%% file: modsurveyfin.tex
\documentclass[10pt,fleqn]{article}
\usepackage{mathptmx}
\usepackage{color}

\input{epsf.tex}
\usepackage{latexsym,amsfonts,amssymb,epsfig,verbatim}
\usepackage{amsmath,amsthm,amssymb,latexsym,graphics,textcomp}
\usepackage{eucal,eufrak}
\usepackage{graphicx}
\usepackage{url}

\theoremstyle{plain}
\newtheorem{theorem}{Theorem}[section]
\newtheorem{lemma}[theorem]{Lemma}
\newtheorem{corollary}[theorem]{Corollary}

\theoremstyle{definition}
\newtheorem*{remark}{Remark}

\newtheorem*{question}{Question}

\newtheorem*{measure}{Measure Conjecture for $\Mod(S)$}
\newtheorem*{gromov}{Gromov's question}

\def\1{\mathbf 1}

\def\GL{\mathcal{GL}}

\def\ML{\mathcal{ML}}

\def\PML{\mathbb P \mathcal{ML}}

\def\WH{\mathfrak{H}}
\def\Mod{\mathrm{Mod}}

\def\Aff{\mathrm{Aff}}

\def\A{\mathfrak A}
\def\C{\mathcal C}
\def\F{\mathcal F}
\def\M{\mathcal M}

\def\T{\mathcal T}

\def\Re{\mathrm{Re}}
\def\Im{\mathrm{Im}}

\def\diam{\mathrm{diam}}

\def\d{\mathrm{d}}

\def\co{\colon\thinspace}

\begin{document}

\title{\textbf{Subgroups of mapping class groups \\ from the geometrical viewpoint}}
\author{Richard P. Kent IV\thanks{The first author was supported by a Donald D. Harrington Dissertation Fellowship.} \  and  Christopher J. Leininger}
\maketitle

\begin{enumerate}
\item[] {\small Once it is possible to translate any particular proof from one theory to another, then the analogy has ceased to be productive for this purpose; it would cease to be at all productive if at one point we had a meaningful and natural way of deriving both theories from a single one. . . . Gone is the analogy: gone are the two theories, their conflicts and their delicious reciprocal reflections, their furtive caresses, their inexplicable quarrels; alas, all is just one theory, whose majestic beauty can no longer excite us.}

\begin{flushright} {\small --- \textit{from Andr\'e Weil's 1940 letter to Simone Weil}}
\footnote{Martin H. Krieger's translation \cite{weilletter} from the French, pages 244--255 of \cite{weilletterfrench}.}
\end{flushright}
\end{enumerate}

\section{The reluctant analogy}

We survey here the analogy between Kleinian groups and subgroups of the mapping class group $\Mod(S)$ of an analytically finite hyperbolic surface $S$. This typically compares the action of a Kleinian group on $\mathbb H^n$ with the action of a subgroup of $\Mod(S)$ on the Teichm\"uller space $\T(S)$ of $S$---itself oft replaced by W. Harvey's complex $\C(S)$ of curves of $S$.  This is often problematic as the Teichm\"uller metric on $\T(S)$ is in dissonance with all reasonable notions of hyperbolicity---though harmonious with Weil's sentiment. The replacement of $\T(S)$ by $\C(S)$ is no sure solution as $\C(S)$, though hyperbolic, is a space nightmarish in complexity---even locally. On occasion, the fact that the points of $\T(S)$ and its boundary are themselves geometric objects allows the transfer of intuition from one theory to another. On rare occasions, this additional level of complexity allows for theorems in $\Mod(S)$ stronger than their echoes in the realm of Kleinian groups.

We make no claim to completeness in the exposition, written as if in conversation with a fellow of the study of Kleinian groups (or that of Teichm\"uller theory)---which, as a consequence, strikes most keenly subjects with which we are most familiar; and, as if at the blackboard, we proceed in relaxed fashion in hopes of transmitting geometric intuition.


\section{The dynamical viewpoint}

Inspired by work of H. Masur \cite{masurhandle}, J. McCarthy and A. Papadopoulos \cite{mccarthypapa} have examined the dynamic behavior of subgroups of $\Mod(S)$ acting on Thurston's space $\PML(S)$ of projectivized measured laminations on $S$.

Let $G$ be a subgroup of $\Mod(S)$. McCarthy and Papadopoulos say that $G$ is \textbf{sufficiently large} if it contains two independent pseudo-Anosov mapping classes---meaning there are two pseudo-Anosov mapping classes whose associated laminations are pairwise distinct.

To simplify the ensuing discussion, we assume throughout that $G$ is sufficiently large, unless otherwise stated---this is most often unnecessary, although the limit set must be suitably interpreted, see \cite{mccarthypapa,KL}.

\subsection{The limit set} \label{limitsetsect}

There is a unique non-empty closed $G$--invariant set $\Lambda_G$ in $\PML(S)$ on which $G$ acts minimally, called \textbf{the limit set}---uniqueness requires the sufficiently large hypothesis.  The limit set is the closure of the set of stable laminations of pseudo-Anosov mapping classes in $G$. It is perfect and has empty interior provided it is not equal to the entire sphere $\PML(S)$---see \cite{mccarthypapa} for proofs of these facts.

\subsubsection*{Example: Affine groups and Veech groups.} Let ${\mathbb H}_q$ be the Teichm\"uller disk associated to a quadratic differential $q$: a totally geodesic hyperbolic plane in $\T(S)$ of constant curvature $-4$.  The visual boundary of ${\mathbb H}_q$ is canonically identified with $\PML(q)$, the projective space of laminations underlying the space of vertical foliations of complex multiples of $q$.  Being identified with the boundary of a hyperbolic plane, $\PML(q)$ is naturally a real projective line.

A class of subgroups of $\Mod(S)$ that have received much attention are stabilizers of $\mathbb H_q$.
Such a group is realized as the subgroup of $\Mod(S)$ which acts by affine homeomorphisms with respect to the flat metric
induced by $q$, and is called an \textbf{affine group} $\Aff(q)$.  If it acts with finite
covolume on $\mathbb H_q$, then $\Aff(q)$ is said to be a \textbf{Veech group}.  For $G < \Aff(q)$, we let $\Lambda_G(q)
\subset \PML(q)$ denote the limit set of $G$ acting as \emph{a Fuchsian group} on $\mathbb H_q \cup \PML(q)$.

\begin{theorem}\label{projective} For any quadratic differential $q$, the obvious inclusion
\[
\PML(q) \to \PML(S)
\]
is a piecewise projective embedding.  In particular, for any $G < \Aff(q)$, the restriction to
$\Lambda_G(q)$ parameterizes the limit set in $\PML(S)$.
\end{theorem}
\begin{proof} We begin by explicitly describing the projective structure on $\PML(q)$.

For simplicity, we suppose that $q$ is the square of an abelian differential (holomorphic $1$--form) $\omega$ and write $\omega = \Re(\omega) + i \Im(\omega)$.  If $\eta$ is any (harmonic) $1$--form in the 2--dimensional span $\langle \Re(\omega), \Im(\omega) \rangle$, then we let $\F_\eta$ denote the measured foliation with foliation tangent to $\ker(\eta)$ and transverse measure given by integrating the absolute value of $\eta$.  The vertical foliations of complex multiples of $q$ are precisely the foliations $\F_\eta$ for $\eta \in \langle \Re(\omega),\Im(\omega) \rangle$.  The map
\[
\langle \Re(\omega),\Im(\omega) \rangle - \{ 0 \} \to \PML(q)
\]
taking $\eta$ to $[\, \F_\eta \,]$, the projective class of $\F_\eta$, has the property that $\eta$ and $\eta'$ lie in the same fiber if and only if $\eta = t \eta'$
for some $t \in \mathbb R$, $t \neq 0$.  This identifies $\PML(q)$ with ${\mathbb P}\langle \Re(\omega),\Im(\omega)
\rangle$, which describes the natural projective structure on $\PML(q)$.

To show that the map from $\PML(q)$ to $\PML(S)$ is piecewise projective, it suffices to show that
\[
\langle \Re(\omega),\Im(\omega) \rangle \to {\mathbb R}^{\mathcal S}
\]
given by $\eta \mapsto \{i(\F_\eta,\gamma)\}_{\gamma \in {\mathcal S}}$ is piecewise linear.  Therefore, we must show $\eta \mapsto i(\F_\eta,\gamma)$ is piecewise linear for any curve $\gamma$.  Representing $\gamma$ by its
$q$--geodesic representative, we see that it is a concatenation $\gamma = \gamma_1 \gamma_2 \cdots \gamma_k$, with each
$\gamma_j$ a $q$--straight segment connecting zeros of $q$.  It follows that
\[
i(\F_\eta,\gamma) = \sum_{j=1}^k |\eta \cdot \gamma_j|
\]
where $\eta \cdot \gamma_j$ is the evaluation of the $1$--form $\eta$ on the segment $\gamma_j$, with respect to either
orientation.  Since $\eta \to \eta \cdot \gamma_j$ is linear in $\eta$, it follows that $i(\F_\eta,\gamma)$ is piecewise
linear in $\eta$, as required.

The general case follows by considering the appropriate double branched covers carefully.
\end{proof}


We note that in general, $\PML(q)$ disagrees with the set of accumulation points of the associated Teichm\"uller disk in Thurston's compactification,
even when $\Aff(q)$ is a Veech group.  This is a consequence of Masur's description in
\cite{twoboundaries} of how Teichm\"uller's visual compactification of $\T(S)$ by $\PML(S)$ differs from Thurston's
compactification.

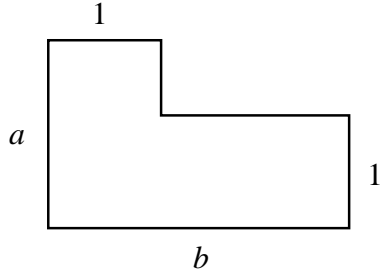
\begin{figure}
\begin{center}
\input{Ltable.pstex_t}
\end{center}
\caption{An $\mathrm L$--shaped table}\label{Ltable}
\end{figure}

As a particular instance, consider an $\mathrm L$--shaped Euclidean polygon $P$ as in figure \ref{Ltable}.  After
dividing each of the sides of lengths $a$ and $b$ appropriately into two pieces, there is an obvious way to identify
parallel sides to obtain a surface $\Sigma$ of genus $2$.  If the gluing maps are simply restrictions of translations,
then $\Sigma$ is naturally a Riemann surface and the differential $z\thinspace \d z$ yields an abelian differential
$\omega$ on $\Sigma$. Squaring $\omega$, we obtain a quadratic differential $q = \omega^2$. C. McMullen and K. Calta
independently proved that the stabilizer of $\mathbb H_q$ will act with finite covolume precisely when $a$ and $b$ are
rational or $a=x + z \sqrt{d}$ and $b = y + z \sqrt{d}$ for some $x,y,z$ in $\mathbb Q$ with $x+y=1$ and $d$ a
nonnegative integer \cite{mcmullenbilliard,calta}. It is easy to see that for most values of $a$ and $b$, the point in
$\PML(q)$ represented by $q$ is a union of simple closed curves with \textit{distinct} weights. On the other hand, the
Teichm\"uller geodesic ray corresponding to $q$ converges in $\T(S) \cup \PML(S)$ to a union of simple closed curves
\textit{all of whose weights are the same}
\cite{twoboundaries}.\\

\noindent Let $\alpha$ and $\beta$ be two simple closed curves that bind $S$. Let $G(\alpha,\beta)$ be the group
generated by Dehn twists $T_\alpha$ and $T_\beta$ in $\alpha$ and $\beta$. One may build from $\alpha \cup \beta$ a
branched flat metric on $S$ ready--equipped with a quadratic differential $q$. Moreover, the Teichm\"uller disk
$\mathbb H_q$ is stabilized by the group $G(\alpha,\beta)$.  We note that in general $G(\alpha,\beta)$ is \textit{not} the full stabilizer of $\mathbb H_q$.  If $i(\alpha,\beta) \geq 3$ (which is true, for instance, when $S$ is closed), then $\Lambda_{G(\alpha,\beta)}(q)$ is a Cantor set in $\PML(q)$; see \cite{leiningertwist}, Section 6.1.  By Theorem \ref{projective} we find
that $\Lambda_{G(\alpha,\beta)}$ is a tame Cantor set.

\subsubsection*{Example: the handlebody groups.}  Let $H$ be a $3$--dimensional handlebody whose boundary is homeomorphic to $S$ and let $\Mod(H)$ be the subgroup of $\Mod(S)$ consisting of mapping classes of $S$ that extend to $H$. H. Masur initiated the study of limit sets of subgroups of $\Mod(S)$ in his study of $\Mod(H)$ \cite{masurhandle}.  In particular, he proved the following theorem.

\begin{theorem}[Masur] The set $\Lambda_{\Mod(H)}$ is connected.\qed
\end{theorem}

\noindent Not only is $\Lambda_{\Mod(H)}$ nowhere dense, but we have the following remarkable theorem of H. Masur \cite{masurhandle} in genus 2 and S. Kerckhoff \cite{kerckhoffhandle} in general.

\begin{theorem}[Masur, Kerckhoff] The set  $\Lambda_{\Mod(H)}$ has measure zero.\qed
\end{theorem}

\bigskip

\noindent The Ahlfors Measure Conjecture \cite{ahlforsmeasure}, which is now known \cite{agol,calgabai,canaryends}, asserts that the limit set of a Kleinian group has zero or full measure in the Riemann sphere.  We have the same conjecture for subgroups of $\Mod(S)$:

\begin{measure} If $G$ is a finitely generated subgroup of $\Mod(S)$, then $\Lambda_G$ has zero or full measure in $\PML(S)$.
\end{measure}
We will see later that this holds for a large class of groups.\\

\noindent Another of the outstanding conjectures in the theory of Kleinian groups is that the limit set of a freely
indecomposable Kleinian group is locally connected.  While there has been a great deal of progress on this conjecture,
it remains open; see
\cite{AMlconnect,ADPlconnect,floydlconnect,minskylconnect,klarreichlconnect,mcmullenlconnect,br1,br2,br}. In the spirit
of curiosity, we ask:

\begin{question}[Local connectivity of $\Lambda$] If $G$ is a one--ended subgroup of $\Mod(S)$, is $\Lambda_G$ locally connected?
\end{question}

\noindent The approach often taken to establish local connectivity for Kleinian groups is to produce a Cannon--Thurston
map from the boundary of the group to the Riemann sphere. Unfortunately, this avenue is closed for us.  Let $\dot S$
denote the result of puncturing the closed surface $S$ once.

\begin{theorem} For any once-punctured surface $\dot S$ there are one--ended hyperbolic subgroups $G$ of $\Mod(\dot S)$ that do not admit Cannon--Thurston maps, meaning that there is no continuous $G$--equivariant extension of $G \to \T(S)$ to $G \cup \partial G \to \T(\dot S) \cup \PML(\dot S)$.
Indeed, for these examples there is no continuous $G$--equivariant map $\partial G \to \PML(\dot S)$.
\end{theorem}

\begin{proof} There is a forgetful map
\[
\Mod(\dot S) \to \Mod(S)
\]
whose kernel is $\pi_1(S)$, see \cite{birman} or section \ref{extensions}.  Since $\pi_1(S)$ is normal, its limit set is all of $\PML(\dot S)$, and so a Cannon--Thurston map $\Phi \co \partial \pi_1(S) \to \PML(\dot S)$ would be surjective.  An element $f$ in $\pi_1(S)$ represented by a simple closed curve $\alpha$ in $S$ fixes pointwise some $\PML(Y)$, where $Y$ is a subsurface of $\dot S$. Let $x$ and $y$ be the fixed points of $f$ in $\partial \pi_1(S)$. Pick a point $\lambda$ in $\PML(Y)$ that misses $\Phi(\{x,y\})$. By equivariance, $f^n \lambda$ must converge to a point in $\Phi(\{x,y\})$ and yet $f^n \lambda = \lambda$, a contradiction.
\end{proof}

It is shown in \cite{graphsofveech} that the closed surface groups $G$ in $\Mod(S)$ obtained from the Combination
Theorem in \cite{LRcombination} (see Section \ref{combination}) also provide examples for which there is no
Cannon--Thurston map $\partial G \to \PML(S)$. However, unlike the example given in the proof, it is also shown that
there is a $G$--equivariant map $\partial G \to \GL(S)$, where $\GL(S)$ is the space of geodesic laminations with the
Thurston topology.

\subsection{The domain of discontinuity}

When accustomed to working with Kleinian groups, one is tempted to consider the action of $G$ on the complement of the
limit set $\PML(S) - \Lambda_G$ in hopes of finding a domain of discontinuity. As Masur observed for the handlebody
group \cite{masurhandle}, $G$ need not act properly discontinuously on $\PML(S) - \Lambda_G$. For example, consider the
groups $G(\alpha,\beta)$ constructed above. The Dehn twist $T_\alpha$ fixes, in particular, every point in $\PML(S -
\alpha)$ and since $\Lambda_G$ is a Cantor set, $T_\alpha$ is an infinite order element fixing a point in $\PML(S) -
\Lambda_G$.

The natural remedy is to enlarge $\Lambda_G$.  Define the \textbf{zero locus $Z \Lambda_G$ of $\Lambda_G$} by
\[
Z \Lambda_G = \{\mu \in \PML(S)\, |\, i(\mu,\lambda) = 0\ \mbox{for\ some\ } \lambda \in \Lambda_G \}.
\]
This was Masur's solution in \cite{masurhandle}.  McCarthy and Papadopoulos prove that for any $G$, $\Delta_G = \PML(S)
- Z \Lambda_G$ is a \textbf{domain of discontinuity} for $G$, meaning that $G$ acts properly discontinuously on
$\Delta_G$ \cite{mccarthypapa}.

\subsubsection*{Example: the $G(\alpha,\beta)$ again.}
We return to the examples $G=G(\alpha,\beta)$ to discuss $Z \Lambda_{G}$.  Given a simple closed curve $\gamma$ in
$\Lambda_G$, not only is the entirety of $\PML(S-\gamma)$ contained in $Z\Lambda_G$, the entire cone
\[
\{ [\mu] \ | \ \mu = t\gamma + \nu \ \mathrm{where}\ t \in [0,\infty) \ \mathrm{and}\ \nu \in \ML(S-\gamma) \}
\]
lives there, too.

In general, given a lamination $\lambda$ in the limit set, we will see the join of the simplex of measures on $\lambda$
and the projective space of measured laminations on its complement.  For the groups $G(\alpha,\beta)$, we will only see
cones, as every limit point of $G$ is uniquely ergodic; see below.  We say that a lamination is \textbf{uniquely
ergodic} if it admits only one transverse measure up to scale, and do not require that it be filling.

For simplicity, let us now assume that $S$ is a five--times punctured sphere.  In this case $\PML(S)$ is the
three--sphere, and we can imagine a caricature of $Z \Lambda_G$. The limit set $\Lambda_G$ is a Cantor set, which we
imagine lying on a great circle in the three--sphere---for reference, we imagine the sphere equipped with its metric of
constant curvature $1$. The rational points $\mathcal X$ of this Cantor set correspond to simple closed curves on $S$,
and at each of them we see a disk (the cone mentioned above). The set of irrational points $\mathcal Y$ consists
entirely of filling uniquely ergodic laminations. To see this, note that any limit point of $G$ is the end of a ray in
the convex hull in $\mathbb H_q$ of $\Lambda_G$. Such a ray either exits a cusp when pushed to $\mathbb H_q/G$, or
returns infinitely often to a fixed compact set. In the first case, the limit point is a conjugate of $\alpha$ or
$\beta$.  In the latter case, the lamination is filling and uniquely ergodic by Masur's criterion
\cite{masurhaus}---compare with the proof of the Veech dichotomy \cite{masurtabachnikov}, Chapter 5.

As a sequence $x_n$ in $\mathcal X$ converges to a point $y$ in $\mathcal Y$, the disks over the $x_n$ shrink and
converge to $y$. To see this, note that each $x_n$ has zero intersection number with every lamination in its cone, and
so any lamination in the latter must converge to a lamination having zero intersection with, and therefore equaling,
$y$.

\subsubsection*{Example: lifting $\Mod(S)$ to covers.}
It is beneficial to pause and contemplate the following example. Let $f \co \widetilde S \to S$ be a finite covering map, and
let $\Mod(f \co \widetilde S \to S)$ be the subgroup of $\Mod(\widetilde S)$ consisting of lifts of mapping classes on $S$.  This is virtually isomorphic to $\Mod(S)$.
The limit set in this case is obtained as the image $f^* \co \PML(S) \to \PML(\widetilde S)$, since this is a closed invariant set containing the fixed points of pseudo-Anosov mapping class as a dense set.  The zero locus is thus a thorny array of joins of spheres and simplices of measures, all of varying sizes and shapes in kaleidoscopic arrangement.\\

\noindent While a group $G$ always acts properly discontinuously on $\Delta_G$, it is not obvious that $\Delta_G$ should be the largest open set for which this is true.
When $Z \Lambda_G = \Lambda_G$, then it is immediate that $\Delta_G$ is a maximal domain of discontinuity.
However, even when $Z \Lambda_G$ properly contains $\Lambda_G$ it can be the case that $\Delta_G$ is maximal.  For example, this is true when $G$ is a Veech group---every point of $Z \Lambda_G - \Lambda_G$ is fixed by some infinite order element of $G$.

\begin{question}[Maximality of $\Delta$] Given a group $G < \Mod(S)$ is $\Delta_G$ the maximal domain of discontinuity for $G$?
\end{question}

The action of $G$ on the preimage of $\Delta_G$ in $\ML$ is also properly discontinuous by general principles.  It is worth pointing out that C. Lecuire has shown that for the handlebody group $\Mod(H)$ the preimage of $\Delta_{\Mod(H)}$ is \textit{not} a maximal domain of discontinuity.\\

Subgroups of $\Mod(S)$ deviate from their cousins in $\mathrm{Isom}(\mathbb H^n)$ in at least one desirable way: the action on the domain of discontinuity puts strong constraints on their geometry.  For instance, a Kleinian group may act cocompactly on its domain of discontinuity and yet still contain parabolic elements, while subgroups of the mapping class group may do no such thing:

\begin{theorem}[\cite{KL}]\label{deltapseudo} If $G < \Mod(S)$ has a non-empty domain of discontinuity $\Delta_G$ on which $G$ acts cocompactly, then $G$ is virtually purely pseudo-Anosov.
\end{theorem}
We will discuss the proof in the following section.


\section{The view from Teichm\"uller space}

\subsection{The behavior of Teichm\"uller geodesics}\label{behaveteichsect}

The \textbf{$\epsilon$--thick part} of Teichm\"uller space is the set of hyperbolic structures on $S$ whose injectivity radius is greater than $\epsilon$.  The \textbf{$\epsilon$--thin part} is the complement of the $\epsilon$--thick part.  A subset of $\T(S)$ is \textbf{$\epsilon$--cobounded} if it lies in the $\epsilon$--thick part, and \textbf{cobounded} if it lies in some $\epsilon$--thick part.

A study of Teichm\"uller geodesics is prerequisite to geometric study of groups acting on the Teichm\"uller space. Especially important is an understanding of the interaction of geodesics with the thick and thin parts of $\T(S)$. K. Rafi's bounded geometry theorem for Teichm\"uller geodesics, in the spirit of Minsky's Bounded Geometry Theorem for geometrically infinite hyperbolic $3$--manifolds, ferries such understanding. To describe this theorem, we will need a host of definitions, which occupy the two sections to follow.

The behavior of Teichm\"uller geodesics is often best understood via the associated laminations, which we now describe.
A pair of binding laminations $\lambda_-,\lambda_+$ in $\PML(S)$ is uniquely associated to a Teichm\"uller geodesic
$\tau(\lambda_-,\lambda_+)$. Conversely, every Teichm\"uller geodesic uniquely determines a pair of laminations
$\lambda_-,\lambda_+$, its \textbf{negative} and \textbf{positive directions}, respectively.  These are the laminations
underlying the horizontal and vertical foliations of any quadratic differential defining the geodesic; see
\cite{hubbardmasur}.

\subsubsection{Complexes of curves and arcs}\label{complexsect}

Let $Y$ be a compact surface.
 Harvey's \textbf{complex of curves} $\C(Y)$ of $Y$ is the simplicial complex whose $k$--cells are collections of isotopy classes of $k+1$ disjoint pairwise non-isotopic essential simple closed curves---when the Euler characteristic of $Y$ is close to zero, other definitions are often required, but we ignore this here for simplicity.

It is convenient to have a \textbf{complex of arcs} $\mathcal A(Y)$.
We define $\mathcal A(Y)$ to be the simplicial complex whose $k$ cells are collections of isotopy classes of $k+1$ disjoint pairwise non-isotopic essential simple closed curves and arcs---where isotopy classes of arcs are defined relative to $\partial Y$.

If $Y$ is an annulus, we need a better definition, and we define $\mathcal A(Y)$ to be the graph whose vertices are isotopy classes of essential arcs  in $Y$ \textit{relative to their endpoints in $\partial Y$} and that two vertices are joined by an edge if they may be realized disjointly.

For any $Y$, we metrize $\C(Y)$ and $\mathcal A(Y)$ by demanding that any simplex is a regular Euclidean simplex with all side lengths equal to one and taking the induced path metric.

Although the complex $\mathcal A(Y)$ is uncountable when $Y$ is an annulus, it is nonetheless quasiisometric to $\mathbb Z$ \cite{MM2}.

\subsubsection{Ivanov--Masur--Minsky subsurface projections}\label{subsurfacesect}

A subsurface $Y$ of $S$ is said to be a \textbf{proper domain} if it is not equal to $S$ and the induced map on fundamental groups is injective.

Fix a hyperbolic metric on $S$ and realize every element of $\PML(S)$ as a geodesic lamination there. Given a proper domain $Y$ in $S$ with geodesic boundary and a geodesic lamination $\lambda$, intersecting $\lambda$ with $Y$ yields a simplex in $\mathcal A(Y)$---when $Y$ is an annulus, a different procedure is required. This simplex is the \textbf{projection of $\lambda$ to $Y$}, denoted $\pi_Y(\lambda)$---note that $\pi_Y(\lambda)$ is allowed to be empty.

Given two geodesic laminations $\mu$ and $\lambda$, the \textbf{projection coefficient for $\mu$ and $\lambda$ at $Y$} is defined to be
\[
\d_Y(\mu,\lambda) = \mathrm{diam}_{\mathcal A(Y)}(\pi_Y(\mu) \cup \pi_Y(\lambda))
\]
\textit{provided that $\pi_Y(\mu) \neq \emptyset$ and $\pi_Y(\lambda)\neq \emptyset$}.  If either of $\pi_Y(\mu)$ or $\pi_Y(\lambda)$ is empty, we define $\d_Y(\mu,\lambda)=\infty$.

We henceforth write $\mathrm{diam}_Y(\, \cdot \,)$ to denote $\mathrm{diam}_{\mathcal A(Y)}(\, \cdot \,)$.



\subsubsection{Bounded geometry theorems}\label{bddgeomsect}

Minsky's Bounded Geometry Theorem \cite{minskybound} says that a geometrically infinite hyperbolic $3$--manifold homeomorphic to $S \times \mathbb R$ has its injectivity radius bounded below if and only if the subsurface projection coefficients of its ending laminations are all uniformly bounded above.

K. Rafi has characterized the short curves in hyperbolic structures on a Teichm\"uller geodesic in terms of intersection data for the subsurface projections of its defining laminations---see \cite{rafi} for the precise statement.  With the global connection between intersection numbers and subsurface projection coefficients described in \cite{MM2}, this yields the following bounded geometry theorem for Teichm\"uller geodesics, which is implicit in the the proof of Theorem 1.5 of \cite{rafi}.

\begin{theorem}[Rafi's Bounded Geometry Theorem] \label{rafitheorem}
For every $D > 0$, there exists $\epsilon >0$ such that if $\tau = \tau(\lambda_-, \lambda_+)$ is a Teichm\"uller geodesic with negative and positive directions $\lambda_-$ and $\lambda_+$ in $\PML(S)$ satisfying
\[
\d_Y(\lambda_-,\lambda_+) \leq D
\]
for every proper domain $Y$ (not a pair of pants) in $S$, then $\tau$ is $\epsilon$--cobounded.

Conversely, for every $\epsilon  > 0$ there exists $D > 0$ such that if $\tau$ is an $\epsilon$--cobounded Teichm\"uller geodesic with negative and positive directions $\lambda_-$ and $\lambda_+$, then
\[
\d_Y(\lambda_-,\lambda_+) \leq D
\]
for every proper domain $Y$ which is not a pair of pants.\qed
\end{theorem}

Another useful theorem is Masur and Minsky's Bounded Geodesic Image Theorem \cite{MM2}:

\begin{theorem}[Masur and Minsky's Bounded Geodesic Image Theorem] \label{mmbgi}
There exists a constant $M = M(S)$ with the following property.
Let $Y$ be a proper domain of $S$ which is not a pair of pants and let $\gamma$ be a geodesic segment, ray, or biinfinite line in $\C(S)$, such that $\pi_{Y}(v) \neq \emptyset$ for every vertex $v$ of $\gamma$.
Then
\[
\diam_{Y}(\gamma) \leq M.
\]\qed
\end{theorem}

\subsection{The hull}

A convenient notion in the theory of Kleinian groups is that of a convex hull. This is the convex hull in $\mathbb H^n$
of a Kleinian group's limit set.  This is always defined, and a Kleinian group is said to be geometrically finite if it
acts with finite covolume on (the $1$--neighborhood of) its convex hull, convex cocompact if it acts cocompactly there.
The situation is not so simple in $\Mod(S)$.  We begin with the notion of a hull.

Let $\A$ be a closed subset of $\PML(S)$.  If $\A$ has the property that for every $\lambda \in \A$, there exists a
$\mu \in \A$ such that $\lambda$ and $\mu$ bind $S$, then we define the \textbf{weak hull} $\WH_\A$ of $\A$ to be the
union of all geodesics in $\T(S)$ whose directions lie in  $\A$.  If $\A$ does not have this property then we say that
the weak hull is not defined.  A set $\WH$ is a \textbf{weak hull} if it is $\WH_\A$ for some closed $\A \subset
\PML(S)$ with the aforementioned property.

\begin{remark} The notion of a weak hull for $G$ in $\Mod(S)$ is due to Farb and Mosher \cite{FMcc}, although we warn the reader that our definition is less restrictive. E. Swenson uses a similar device in his study of groups of isometries of Gromov hyperbolic metric spaces \cite{swenson}.
\end{remark}

Since we are assuming $G$ sufficiently large, $\Lambda_G$ has a weak hull $\WH_G$: $\Lambda_G$ contains the stable
lamination $\mu$ of a pseudo-Anosov mapping class, which will bind with any other lamination in $\Lambda_G$.

A leitmotif of Teichm\"uller theory is \textit{the thick part feigns thinness (hyperbolicity)}. An instance:

\begin{theorem}[\cite{KL}]\label{thintriangles} A triangle in Teichm\"uller space with cobounded sides is thin. \qed
\end{theorem}
\noindent The proof analyzes triangles directly,  applications of Minsky's Quasiprojection Theorem \cite{minskycrelle}
and Masur's Asymptotic Rays Theorem \cite{masuruniquely}  providing thinness.  M. Duchin has discovered a different
proof \cite{duchin}.
The following theorem is an easy application of Theorem \ref{thintriangles}.

\begin{theorem}[\cite{KL}]\label{hyperbolichull} Let $\WH=\WH_\A$ be a cobounded weak hull.  Then $\WH$ is quasiconvex and there is a metric neighborhood of $\WH$ which is a $\delta$--hyperbolic metric space---when equipped with its induced path metric. \qed
\end{theorem}

\begin{theorem}[\cite{KL}]\label{coboundedhull}
If $G < \Mod(S)$ has a non-empty domain of discontinuity $\Delta_G$ on which $G$ acts cocompactly, then $\WH_G$ is cobounded and so every lamination in $\Lambda_G$ is filling and uniquely ergodic.
\end{theorem}
\begin{proof}[Sketch] The proof begins with a demonstration that every lamination in $\Lambda_G$ is filling, and culminates with a generalization of Minsky's proof when $G$ is cyclic \cite{minskyinj}.

Define
\[
Z Z \Lambda_G = \{\mu \in \PML(S)\, |\, i(\mu,\lambda) = 0\ \mbox{for\ some\ } \lambda \in Z\Lambda_G \}.
\]

One first shows that $Z  Z\Lambda_G = Z \Lambda_G$.  For this, we find a simple closed curve $\alpha \in Z Z \Lambda_G - Z \Lambda_G$ (assuming this to be nonempty) having zero intersection number with some $\mu$ in $Z \Lambda_G$.  The $1$-simplex of measures defined by $\{\nu_t = t \alpha + (1-t)\mu \, | \, t \in [0,1]\}$ meets $\Delta_G$ in the set $\{ \nu_t \, | \, t \in (0,1] \}$.  Moreover, this noncompact segment is properly embedded in $\Delta_G$, hence closed.  A compact fundamental domain for $\Delta_G$ would require infinitely many translates to cover this set, and the particular construction of a fundamental domain we use (see \cite{mccarthypapa}) shows that this is impossible.

One then proves that every lamination in $\Lambda_G$ is filling as follows.
If there were a non-filling lamination in $\Lambda_G$ there would be a simple closed curve $\alpha$ in $Z \Lambda_G$. Now, since $Z \Lambda_G = Z Z \Lambda_G$, any closed curve disjoint from $\alpha$ is in $Z \Lambda_G = Z Z \Lambda_G$, and we conclude that every simple closed curve lies in $Z \Lambda_G$, as the complex of curves is connected.  As these are dense in $\PML(S)$, we conclude that $Z\Lambda_G = \PML(S)$, contradicting the fact that $\Delta_G$ is non-empty.

Given a compact set $\mathcal K$ in $\Delta_G$, $\kappa$ in $\mathcal K$, and $\lambda$ in $\Lambda_G$, the
intersection $\lambda \cap (S-\kappa)$ is a collection of arcs, as every lamination in $\Lambda_G$ is filling. The
hyperbolic length of the arcs of $\lambda \cap (S-\kappa)$ may be bounded over all $\lambda$ in $\Lambda_G$ and
$\kappa$ in $\mathcal K$, as $\Lambda_G$ and $\mathcal K$ are compact.  If $\mathcal K$ is a fundamental domain for the
action of $G$ on $\Delta_G$, equivariance provides a bound on the hyperbolic length of $\lambda \cap (S - \kappa)$ for
all $\lambda$ in $\Lambda_G$ and $\kappa$ in $\Delta_G$.  In particular, for any subsurface $Y$ (realized as the
interior of a surface with geodesic boundary) and any $\lambda \in \Lambda_G$ there is a uniform bound on the length of
any segment of $Y \cap \lambda$.  The Collar Lemma then provides a bound on all of the Masur--Minsky subsurface
projection coefficients $\mathrm d_Y(\lambda_-,\lambda_+)$, where $\lambda_-$ and $\lambda_+$ are in $\Lambda_G$---when
$Y$ is an annulus, an alternate argument is required. With Rafi's bounded geometry theorem, this allows us to conclude
that geodesics whose directions lie in $\Lambda_G$ are uniformly cobounded, and hence filling and uniquely ergodic by
Masur's criterion \cite{masurhaus}.
\end{proof}

Theorem \ref{deltapseudo} follows quickly from this theorem.
\begin{proof}[Proof of Theorem \ref{deltapseudo}] Since $\Mod(S)$ possesses a torsion free subgroup of finite index, so does $G$, and so, to complete the proof, it suffices to show that $G$ contains no reducible element. But iterating such an element on $\Lambda_G$ would produce a non-filling lamination in $\Lambda_G$, which is excluded by the conclusion of Theorem \ref{coboundedhull}.
\end{proof}

\subsection{The Measure Conjecture}\label{measureconjsection} Here we verify the measure conjecture for a certain class of subgroups of
$\Mod(S)$.  Indeed, this is a corollary of the work of Masur \cite{masurinterval}, as we now explain.

We say that a projective measured lamination $\lambda$ is $\epsilon$--cobounded if every Teichm\"uller geodesic ray
with direction $\lambda$ eventually lies entirely in the $\epsilon$--thick part of Teichm\"uller space.

\begin{theorem}[Masur] The set $\Lambda_\epsilon$ of $\epsilon$--cobounded laminations has measure zero in $\PML(S)$.
\end{theorem}
\begin{proof}  Since $\Mod(S)$ acts by pulling back hyperbolic structures, $\Lambda_\epsilon$ is $\Mod(S)$--invari-ant.  Since the action of $\Mod(S)$ on $\PML(S)$ is ergodic \cite{masurinterval}, $\Lambda_\epsilon$ has zero or full measure.

The unit cotangent bundle of Teichm\"uller space is the bundle of unit norm quadratic differentials $\mathcal Q^1$, and
there is a map from $\Pi \co \mathcal Q^1 \to \PML(S)$, obtained by sending a quadratic differential to the projective
class of its vertical foliation (and passing to the underlying lamination). Moreover, in \cite{masurinterval}, Masur
constructs a $\Mod(S)$--invariant measure $\mathrm m$ on $\mathcal Q^1$ with the following properties: $\mathrm m$
takes positive values on open sets; for any subset $E$ of $\PML(S)$, the measure $\mathrm m(\Pi^{-1}(E))$ is zero if
and only if $E$ has Lebesgue measure zero; the $\mathrm m$--measure of $\mathcal Q^1/\Mod(S)$ is finite; and the
geodesic flow $\varphi_{\thinspace t}$ on $\mathcal Q^1/\Mod(S)$ preserves $\mathrm m$.

Suppose that $\Lambda_\epsilon$ has full measure in $\PML(S)$.  Pulling $\Lambda_\epsilon$ back to $\mathcal Q^1$, we
discover that in almost every direction, the geodesic flow eventually carries us over the $2\epsilon$--thick part for
all time---meaning that, for almost every $q$ in $\mathcal Q^1/\Mod(S)$, we have $\varphi_{\thinspace t}(q)$ projecting
to the $2\epsilon$--thick part of $\M(S)$ for all $t$ above some threshold depending on $q$.  This violates Poincar\'e
Recurrence.  To see this, note that the projection $\mathcal Q^1/\Mod(S) \to \M(S)$ is continuous, so that the inverse
image $\mathcal Q^1_{\epsilon/2}$ of the $\epsilon/2$--thin part is open and hence has positive $\mathrm m$--measure.
By the Poincar\'e Recurrence Theorem, for almost every $q$ in $\mathcal Q^1_{\epsilon/2}/ \Mod(S)$,
$\varphi_{\thinspace t}(q)$ must return to $\mathcal Q^1_{\epsilon/2}/\Mod(S)$ for arbitrarily large values of $t$. But
the $2\epsilon$--thick part and the $\epsilon/2$--thin part do not intersect.
\end{proof}

\begin{corollary}[Measure Conjecture when $\Delta_G/G$ is compact]\label{measurecor} If $G < \Mod(S)$ has a non-empty domain of discontinuity $\Delta_G$ on which $G$ acts cocompactly, then $\Lambda_G$ has measure zero.
\end{corollary}
\begin{proof} By Theorem \ref{coboundedhull}, the limit set $\Lambda_G$
is contained in $\Lambda_\epsilon$ for some positive $\epsilon$.
\end{proof}

\section{Convex cocompactness}

A Kleinian group is convex cocompact if it acts cocompactly on the convex hull in $\mathbb H^3$ of its limit set. B. Farb and L. Mosher have extended this notion to subgroups of $\Mod(S)$ \cite{FMcc}. A finitely generated subgroup of $\Mod(S)$ is said to be \textbf{convex cocompact} if it satisfies one of the conditions in the following theorem.

\begin{theorem} [Farb--Mosher \cite{FMcc}] \label{famoconv}
Given a finitely generated subgroup $G$ of $\Mod(S)$, the following statements are equivalent:
\begin{itemize}
\item Some orbit of $G$ is quasiconvex in $\T(S)$.
\item Every orbit of $G$ is quasiconvex in $\T(S)$.
\item \textbf{(Convex cocompact)} $G$ is
hyperbolic and there is a $G$--equivariant embedding $\partial f\co \partial G \to \PML(S)$ with image $\Lambda_G$ such that the weak hull $\WH_G$ of $\Lambda_G$ is defined; the action of $G$ on $\WH_G$ is cocompact; and, if $f\co G \rightarrow \WH_G$ is any $G$--equivariant map, then $f$ is a quasiisometry and the following map is continuous:
\[
\overline{f} = f \cup \partial f\co G \cup \partial G \to \T(S) \cup \PML(S). \qed \]
\end{itemize}
\end{theorem}

\subsection{Characterization}

As the Teichm\"uller space fails to be hyperbolic in any reasonable sense of the word \cite{masurclass,masurwolf,brockfarb}, and in fact exhibits behavior characteristic of positive curvature \cite{minskyextremal}, it is perhaps somewhat surprising that many of the characterizations of convex cocompact Kleinian groups have analogs in $\Mod(S)$.

\subsubsection{From Teichm\"uller space}

Farb and Mosher's definition of convex cocompactness may be streamlined to more closely resemble the definition for Kleinian groups:

\begin{theorem}[\cite{KL}]\label{cocompactonhull} A finitely generated subgroup $G$ of $\Mod(S)$ is convex cocompact if and only if the weak hull $\WH_G$ is defined and $G$ acts cocompactly on $\WH_G$.
\end{theorem}
\begin{proof} One direction is transparent from Theorem \ref{famoconv}.

If $G$ acts cocompactly on $\WH_G$, then the orbit of any point in $\WH_G$ is $A$--dense in $\WH_G$ for some $A$.  Moreover, $\WH_G$ is cobounded.  By Theorem \ref{hyperbolichull}, $\WH_G$ is quasiconvex, and hence so is the orbit of any point in $\WH_G$.
\end{proof}

\subsubsection{From $\PML(S)$}

A group action on a compactum $X$ is a \textbf{discrete convergence action} if the group acts properly discontinuously on the
space of distinct triples $\{(a,b,c) \in X^3\ | \ a \neq b \neq c \neq a \}$ and \textbf{uniform} if the action is
cocompact.  We call a group \textbf{elementary} if it is virtually abelian.

B. Bowditch has characterized hyperbolic groups topologically \cite{bowditchtop}:

\begin{theorem}[Bowditch] A nonelementary group acts as a discrete uniform convergence group on a perfect metrizable compactum $X$ if and only if the group is hyperbolic.  Moreover, $X$ is its Gromov boundary.\qed
\end{theorem}

A Kleinian group is convex cocompact if and only if it acts as a discrete uniform convergence group on its limit set.  This fails in $\Mod(S)$:

\begin{theorem}[\cite{KLtrip}] There are nonabelian free groups in $\Mod(S)$ whose limit sets parameterize their Gromov boundaries (and so act as discrete uniform convergence groups there) and yet fail to be convex cocompact.
\end{theorem}
\begin{proof}[The examples] Let $f$ be pseudo-Anosov and $h$ reducible and pseudo-Anosov on the complement of a nonseparating curve.  After replacing $f$ and $h$ by appropriate powers of conjugates, one constructs an orbit map to $\C(S)$ for which almost every geodesic ray in $\langle f,h \rangle$, namely those that do not end in an infinite sequence of $h$--edges, is sent to a (non-uniform) quasigeodesic ray in $\C(S)$---this construction uses the Bounded Geodesic Image Theorem.

This associates an unmeasured lamination to almost every point of the Gromov boundary $\partial \langle f,h \rangle$.
This is promoted to an equivariant map $\partial \langle f,h \rangle \to \PML(S)$ by appealing to work of Masur
\cite{masurhaus} and Rafi \cite{rafi} to lift the map already defined, and then by an explicit extension to the
remaining set of measure zero.  Continuity follows from various intersection number arguments.

Note that these examples are not convex cocompact since such groups are always purely pseudo-Anosov; see Proposition
3.1 of \cite{FMcc}.
\end{proof}

As is often the case in $\Mod(S)$, the correct replacement for the limit set is its zero locus:

\begin{theorem}[\cite{KLtrip}]\label{convergence} A nonelementary finitely generated subgroup $G$ of $\Mod(S)$ is convex cocompact if and only if $G$ acts as a discrete uniform convergence group on $Z \Lambda_G$.
\end{theorem}
\begin{proof}[Sketch] One first shows that every lamination in $\Lambda_G$ is filling and uniquely ergodic.  Suppose to the contrary that this is not the case. Then there is a positive dimensional simplex $\sigma$ in $Z \Lambda_G$ and one shows that for \textit{any} sequence $g_n$ in $G$, one may pass to a subsequence so that the restriction of $g_n$ to a nontrivial projective linear arc in the interior of $\sigma$ converges uniformly to a (possibly degenerate) arc in $Z\Lambda_G$.  We note that the requirement that the arc lie in the interior of $\sigma$ is necessary to find such a subsequence.
Since the action on $Z \Lambda_G$ is a discrete uniform convergence action, for any $x$ in $Z\Lambda_G$, there is a
sequence $g_n$ so that $g_n(x)$ converges to a point $a$ and $g_n$ restricted to $Z \Lambda_G - \{ x \}$ converges
uniformly on compact sets to a constant function equal to $b \neq a$. If $x$ and $y$ are endpoints of an arc in
$\sigma$ as above, then the limiting arc is non-degenerate, which contradicts uniform convergence to $b$ on $Z
\Lambda_G - \{ x \}$.

A uniquely ergodic filling laminations $\mu$ has the property that if $i(\nu,\mu) = 0$ then $\nu = \mu$.  So as all laminations in $\Lambda_G$ are uniquely ergodic and filling, we note that $\Lambda_G = Z \Lambda_G$.

Because any two distinct uniquely ergodic filling laminations bind $S$ and hence define a Teichm\"uller geodesic, one can construct a continuous equivariant map from the space of distinct triples in $Z \Lambda_G$ to $\T(S)$ by sending $(\lambda,\mu,\nu)$ to the balance point of $\nu$ on the geodesic $\tau(\lambda,\mu)$---this is the unique point where the intersection number of $\nu$ with the vertical foliation is equal to that of the horizontal foliation.  Cocompactness allows one to prove that every point in $\Lambda_G$ is a {\em conical limit point} in the sense described below, and Theorem 4.6 implies $G$ is convex cocompact.
\end{proof}

\subsubsection{Putting the two together}

A convenient notion from the theory of Kleinian groups is that of a conical limit point, which characterizes points in
the limit set by the manner in which they are approached in hyperbolic space.  In our setting: a limit point $\lambda$
in $\Lambda_G$ is said to be \textbf{conical} if every geodesic ray in $\T(S)$ with direction $\lambda$ has a metric
neighborhood containing infinitely many points in a $G$--orbit---we stress to the reader that the definition of conical
is adapted from the geometric definition in Kleinian groups \cite{BeMa}, rather than from the theory of convergence
groups \cite{bowditchtop}.

\begin{theorem}[\cite{KL}]\label{conical} A finitely generated subgroup $G$ of $\Mod(S)$ is convex cocompact if and only if every limit point of $G$ is conical.\qed
\end{theorem}

That all limit points of a convex cocompact group are conical follows easily from the fact that every limit point is
the endpoint of a geodesic in $\WH_G$, and that $G$ acts cocompactly on $\WH_G$ by Theorem \ref{famoconv}.

Generalizing work of McCarthy and Papadopoulos, a fundamental domain for the action of $G$ on $\T(S) \cup \Delta_G$ is constructed in \cite{KL} via F. Bonahon's theory of geodesic currents.  Having only conical limit points implies that this domain is compact by an extremal length argument and so $G$ acts cocompactly on $\T(S) \cup \Delta_G$.  Just as in the hyperbolic setting, this is equivalent to convex cocompactness:

\begin{theorem}[\cite{KL}] A finitely generated subgroup $G$ of $\Mod(S)$ is convex cocompact if and only if $G$ acts cocompactly on $\T(S) \cup \Delta_G$.
\end{theorem}
\begin{proof}[Sketch.] We sketch the proof that if $(\T(S) \cup \Delta_G)/G$ is compact, then $G$ is convex cocompact.  We first observe that the domain $\Delta_G$ is non-empty, as no subgroup of $\Mod(S)$ acts cocompactly on $\T(S)$.  Furthermore, since $\Delta_G/G \subset (\T(S) \cup \Delta_G)/G$ is closed, it is compact, so by theorem \ref{coboundedhull}, $\WH_G$ is cobounded.  Since $(\T(S) \cup \Delta_G ) / G$ is compact, and $\WH_G/G$ is closed therein (this requires some work), $\WH_G/G$ is also compact.  Therefore, $G$ is convex cocompact by Theorem \ref{cocompactonhull}.
\end{proof}

With Corollary \ref{measurecor}, this implies that convex cocompact groups satisfy the Measure Conjecture.

\begin{corollary} If $G$ is convex cocompact, then $\Lambda_G$ has measure zero. \qed
\end{corollary}

\subsubsection{From the complex of curves}


A celebrated theorem of Masur and Minsky \cite{MM1} says that $\C(S)$ is a $\delta$--hyperbolic metric space. Although $\C(S)$ fails to be locally compact, its hyperbolicity often makes it a desirable substitute for the Teichm\"uller space. For instance, a quasiisometric embedding $G \to \T(S)$ is insufficient to guarantee convex cocompactness: a partial pseudo-Anosov mapping class generates a cyclic group quasiisometrically embedded in $\T(S)$ \cite{MMunstable} that fails to be convex cocompact---the orbit is an unstable quasigeodesic; In the complex of curves, quasiisometric embeddings suffice, as proven independently by the authors and U. Hamenst\"adt:

\begin{theorem}[Hamenst\"adt \cite{hamenstadt},  Kent--Leininger \cite{KL}]\label{complextheorem} A finitely generated subgroup $G$ of $\Mod(S)$ is convex cocompact if and only if sending $G$ to an orbit in the complex of curves defines a quasiisometric embedding $\Phi\co G \to \C(S)$.
\end{theorem}
\noindent We sketch the proof from \cite{KL}.
\begin{proof}[Sketch] We first assume that $\Phi \co G \to \C(S)$ is a quasiisometric embedding.  It follows that $G$ is Gromov hyperbolic. Given distinct points $m_-,m_+$ in the Gromov boundary $\partial G$, there is a geodesic in $G$ joining them. This geodesic is carried to a quasigeodesic in $\C(S)$, which is uniformly close to a geodesic $\gamma$ joining $\Phi(m_-)$ and $\Phi(m_+)$.

The proof begins by uniformly bounding the Masur--Minsky subsurface projection coefficients $\d_Y(\Phi(m_-),
\Phi(m_+))$.  To do so it suffices to bound the diameter $\diam_{Y}(\gamma)$ of the projection of $\gamma$ to an
arbitrary subsurface $Y$.

If $Y$ is a proper domain whose boundary is far from $\gamma$, Masur and Minsky's Bounded Geodesic Image Theorem
\cite{MM2} (Theorem \ref{mmbgi} here) provides a bound on $\diam_{Y}(\gamma)$.  If $\partial Y$ is close to $\gamma$, it is close to $\Phi(G)$. In
fact, by translating, we may assume that $\gamma$ and $\partial Y$ are both uniformly close to $\Phi(\1)$. Since the
two ends of $\gamma$ diverge, $\gamma$ may be decomposed into three parts: a finite segment $\gamma_0$ near $\partial
Y$ and two infinite rays $\gamma_-$ and $\gamma_+$ far from $\partial Y$.  The Bounded Geodesic Image Theorem again
bounds $\diam_{Y}(\gamma_\pm)$. The segment $\gamma_0$ fellow travels the image of a geodesic segment in $G$ lying in a
fixed neighborhood of $\1$. Finiteness of this neighborhood allows us to bound $\diam_Y(\gamma_0)$.  The triangle
inequality provides the bound on $\diam_Y(\gamma)$.

Rafi's bounded geometry theorem now shows that the weak hull $\WH_G$ is cobounded.  By Theorem \ref{hyperbolichull},
$\WH_G$ is essentially a hyperbolic metric space, and lifting our $G$--orbit from $\C(S)$ to $\WH_G$, we have a
quasiisometric embedding $G \to \WH_G$.  By the stability of quasigeodesics in hyperbolic metric spaces, we see that
$G$ is quasiconvex in $\WH_G$, and hence in $\T(S)$, since $\WH_G$ is quasiconvex in $\T(S)$ (by Theorem
\ref{hyperbolichull}.)

Now suppose that $G$ is convex cocompact, so that the weak hull $\WH_G$ lies in a thick part of Teichm\"uller space.  A generalization \cite{KL} of Minsky's work on quasiprojections to Teichm\"uller geodesics \cite{minskycrelle} provides a closest points projection $\T(S) \to \WH_G$, and demonstrates that the region of $\T(S)$ where a simple closed curve is short is sent to a uniformly bounded diameter set in $\WH_G$. This projection may be extended to a Lipschitz projection from the electric Teichm\"uller space (quasiisometric to $\C(S)$) to an electrified hull (quasiisometric to $\WH_G$ as we only electrify uniformly bounded diameter sets). We thus obtain a quasiisometric embedding $G \to \C(S)$.
\end{proof}

\begin{corollary} If $S$ is closed, $G$ is a convex cocompact subgroup of $\Mod(S)$, and $G$ lifts to $\Mod(\dot S)$, then its lift is convex cocompact in $\Mod(\dot S)$.
\end{corollary}
\begin{proof} Since $S$ is closed,  the forgetful map $\dot S \to S$ induces a map $\C(\dot S) \to \C(S)$, which is clearly $1$--Lipschitz.  Lifting $G$--orbits from $\C(S)$ to $\C(\dot S)$, we see that  any $G$--orbit in $\C(\dot S)$ is a quasiisometric embedding.
\end{proof}

\subsubsection*{The handlebody group again}

Work of Masur and Minsky \cite{MMdiskset} gives us some geometric information about the mapping class group $\Mod(H)$ of  a handlebody $H$. Letting $\mathcal D$ denote the set of essential simple closed curves in $S=\partial H$ that bound disks in $H$, they prove the following theorem.

\begin{theorem}[Masur--Minsky \cite{MMdiskset}] The set $\mathcal D$ is quasiconvex in the complex of curves $\C(S)$.\qed
\end{theorem}

Although $\Mod(H)$ admits no quasiisometric embedding into $\C(S)$, we do have the following corollary.

\begin{corollary} The handlebody group has a quasiconvex orbit in $\C(S)$. \qed
\end{corollary}

Appealing to the work of Masur and Minsky \cite{MM2}, the second author in \cite{graphsofveech} has obtained similar
quasiconvexity results for the groups from Theorem \ref{LRcombothrm}.

\subsection{Geometrical finiteness?} It is natural to wonder if there is a reasonable notion of geometrical finiteness in $\Mod(S)$.  One's intuition suggests that any definition should encompass both finitely generated affine groups (see section \ref{limitsetsect}) and the combinations of affine groups constructed by the second author and Reid \cite{LRcombination} (see section \ref{combination} following).  We refer the reader to Mosher's article \cite{moshersurvey} for further
discussion and \cite{graphsofveech} for related results.

\section{Kleinian constructions}

\subsection{The Tits alternative and the Schottky argument}

McCarthy has proven that $\Mod(S)$ satisfies the Tits Alternative \cite{mctits}, and the proof follows the classical
Schottky argument on hyperbolic space, showing that two elements have powers that commute or else have powers that
generate a free group. Just as in the Kleinian group setting, this free group may be taken to be convex cocompact:

\begin{theorem}[Mosher \cite{hypbyhyp}, Farb--Mosher \cite{FMcc}] Given two independent pseudo-Anosov mapping classes $\varphi, \psi$, there is a number $\ell$ such that for all natural numbers $m > \ell$, the group generated by $\varphi^m, \psi^m$ is free and convex cocompact.
\end{theorem}
A free convex cocompact subgroup of $\Mod(S)$ is called \textbf{Schottky}.  We discuss two proofs:
\begin{proof}[Sketch of proof from \cite{FMcc}] The proof of McCarthy's Tits Alternative implies that for sufficiently large $m$, the group $\langle \varphi^m, \psi^m \rangle$ is free.
Mosher proved directly in \cite{hypbyhyp} that there is an $\ell$ such that for all $m > \ell$, the canonical $\pi_1(S)$--extension of $\langle \varphi^m, \psi^m \rangle$ (see section \ref{extensions}) is hyperbolic, using the Bestvina--Feighn Combination Theorem.  Theorem \ref{hyperbolicextension} below proves the theorem.
\end{proof}
\begin{proof}[Sketch of proof from \cite{KL}] Given two isometries of a Gromov hyperbolic space acting hyperbolically with pairwise distinct fixed points in the Gromov boundary, one may always raise them to powers so that the resulting group is free and quasiisometrically embedded by its orbits.
Since $\C(S)$ is hyperbolic, an application of Theorem \ref{complextheorem} proves the theorem.
\end{proof}

\subsection{Combination theorems}\label{combination}

The second author and Reid have proven an analog of the first Klein--Maskit combination theorem for affine subgroups of
$\Mod(S)$:

\begin{theorem}[Leininger--Reid \cite{LRcombination}]\label{LRcombothrm} Given two affine groups $H$ and $K$ with a common maximal parabolic subgroup $A$ centralized by a sufficiently complicated mapping class $f$, the group generated by $H$ and $f^{-1}Kf$ is isomorphic to the amalgamated product $H*_A K$.  Moreover, any infinite order element of this group is either pseudo-Anosov or conjugate to a parabolic element in one of the factors. \qed
\end{theorem}

When $S$ is closed, the existence of Veech groups with a single cusp provides the following corollary.

\begin{theorem}[Leininger--Reid \cite{LRcombination}]\label{LRsurfacethrm} If $S$ is closed, then $\Mod(S)$ contains subgroups isomorphic to the fundamental group of a closed surface for which all but one conjugacy class of elements (up to powers) is pseudo-Anosov. \qed
\end{theorem}

In \cite{graphsofveech}, the second author proves a more general combination theorem for certain graphs of Veech groups
and pursues a geometric investigation of their subgroups.  One corollary is the following theorem.

\begin{theorem}[Leininger \cite{graphsofveech}] \label{graphsleininger} Finitely generated virtually purely pseudo-Anosov subgroups of the Leininger--Reid combinations of Veech groups are convex cocompact. \qed
\end{theorem}

A theorem of H. Min provides another geometric combination theorem:

\begin{theorem}[H. Min \cite{min}]\label{min} Given two finite subgroups $H$ and $K$ of $\Mod(S)$, one may conjugate one of the two to obtain a convex cocompact subgroup isomorphic to $H*K$. \qed
\end{theorem}

All the examples of virtually purely pseudo-Anosov subgroups produced by Theorems \ref{graphsleininger} and \ref{min} are virtually free.
This is the case for all known examples.

\section{Global theorems}

\subsection{Abelian subgroups and geometric rank}

The \textbf{geometric rank} of a metric space is the maximal rank of a quasiisometrically embedded Euclidean space. The geometric rank of a Kleinian group is the maximal rank of its free abelian subgroups. J. Behrstock and Minsky \cite{behrstockminsky}, and, independently Hamenst\"adt \cite{hamenstadtrank}, have proven the $\Mod(S)$--analog, answering a question of J. Brock and Farb \cite{brockfarb}:

\begin{theorem}[Behrstock--Minsky \cite{behrstockminsky}, Hamenst\"adt \cite{hamenstadtrank}] The geometric rank of $\Mod(S)$ is the maximal rank of its free abelian subgroups. \qed
\end{theorem}

Both proofs appeal to a study of the asymptotic cone of $\Mod(S)$.

\subsection{Relative hyperbolicity}

A Kleinian group is strongly relatively hyperbolic, in the sense of Farb \cite{relativefarb}, with respect to its maximally parabolic subgroups.  Although we often compare $\Mod(S)$ with a finite covolume Kleinian group, the former does not enjoy the strong relative hyperbolicity of the latter, as established independently by J. Anderson, J. Aramayona, and K. Shackleton \cite{aas}; Behrstock, C. Drutu, and Mosher \cite{bdm}; Bowditch \cite{bowditchcomplexofcurves}; and A. Karlsson and G. Noskov \cite{kn}:

\begin{theorem}[\cite{aas,bdm,bowditchcomplexofcurves,kn}] The mapping class group $\Mod(S)$ is not strongly relatively hyperbolic with respect to any collection of subgroups.\qed
\end{theorem}

\begin{remark} Although Masur and Minsky \cite{MM1} have proven that the Teichm\"uller space is hyperbolic relative to its thin parts, and so $\Mod(S)$ is hyperbolic relative to its maximal abelian subgroups, the ``electrified'' regions fail to satisfy the bounded penetration hypothesis required by strong relative hyperbolicity.
\end{remark}






\section{Extensions}\label{extensions}

We assume from here on that $S$ is closed.

There is a universal extension
\begin{equation} \label{ses1}
1 \to \pi_1(S) \to \Mod(\dot S) \to \Mod(S) \to 1
\end{equation}
corresponding to the universal curve over the moduli space
\begin{equation} \label{quo1}
S \to \M(\dot S) \to \M(S)
\end{equation}
that places above a generic point the Riemann surface it represents (unless $S$ has genus two, where the generic fiber
is the quotient surface by its hyperelliptic involution). Whether or not these extensions are hyperbolic depends on
convex cocompactness of the quotient:

\begin{theorem}[Farb--Mosher \cite{FMcc}]\label{hyperbolicextension} When $F$ is finitely generated and free, $\Gamma_{\! F}$ is hyperbolic if and only if $F$ is Schottky.
\end{theorem}
\begin{proof}[Sketch] Suppose that $F$ is free.  Then the extension is split, and $\Gamma_{\! F}$ is a graph of surface groups. There is an obvious $\mathrm K(\Gamma_{\! F},1)$ homeomorphic to a surface bundle over a finite graph, which we equip with some metric. The Bestvina--Feighn Combination Theorem \cite{combinationtheorem} says, roughly speaking, that $\Gamma_{\! F}$ is hyperbolic if $\Gamma_{\! F}$ has \textit{uniform flaring}, which means that if we push an element of the fundamental group of the fiber around a biinfinite path in the base, its length, in one direction or the other, eventually grows exponentially, or ``flares.''

If $F$ is Schottky, then its hull is cobounded, and given a path in $\WH_F$, this implies that the canonical bundle
over this path possesses the flaring condition in manner depending only on the coboundedness constant. Now, the
universal covers of the canonical bundle over $\WH_G$ and the $\mathrm K(\Gamma_{\! F},1)$ are quasiisometric, and the
Combination Theorem completes the proof.

The other direction is established via S. Gersten's converse \cite{gersten} to the Combination Theorem and the Quasigeodesic Stability Theorem of L. Mosher \cite{mosherstable} and B. Bowditch \cite{bowditchstacks}, which characterizes quasigeodesics close to cobounded geodesics via flaring.
\end{proof}

It is a theorem of Mosher that if $\Gamma_G$ is Gromov hyperbolic, then so is $G$ \cite{mosherextension}, and, after proving a higher dimensional version of Gersten's converse to the Combination Theorem, Farb and Mosher prove that $G$ is in fact convex cocompact \cite{FMcc}:

\begin{theorem}[Farb--Mosher] If the extension $\Gamma_G$ is hyperbolic, then $G$ is convex cocompact. \qed
\end{theorem}

The Bestvina--Feighn Combination Theorem gives sufficient conditions by which a graph of hyperbolic groups is itself hyperbolic.  A higher dimensional version of this theorem (for complexes of groups) is unfortunately not available.  Nonetheless, Hamenst\"adt has proven the following.

\begin{theorem}[Hamenst\"adt \cite{hamenstadt}]\label{ham} If $G$ is a finitely generated subgroup of $\Mod(S)$, the extension $\Gamma_G$ is hyperbolic if and only if $G$ is convex cocompact.\qed
\end{theorem}
\begin{proof}[Sketch]  The key ingredient is a result in the spirit of the Bestvina--Feighn Combination Theorem for certain types of metric fiber bundles over arbitrary hyperbolic spaces, where the fibers are trees.  Namely, a uniform flaring criterion for vertical distances implies hyperbolicity.

Given this, let $\widetilde S$ denote the universal cover of $S$.  One considers an $\widetilde S$--bundle $Y$ over the
Cayley graph ${\frak C}(G)$ of $G$, on which $\Gamma_G$ acts cocompactly.  The Gromov boundary of each fiber is a
circle, and via the action by $\pi_1(S) < \Gamma_G$ on $Y$, one obtains a canonical homeomorphism between the circle of
any given fiber, and the circle $S^1_\infty$ for a particular one.  Now, given a pair of points $x,y \in S^1_\infty$,
from these homeomorphisms we obtain a pair of points in each circle over a point in ${\frak C}(G)$.  The geodesics in
each fiber connecting the pairs of points forms a line sub-bundle of $Y$ over ${\frak C}(G)$.  As in the proof of
Theorem \ref{hyperbolicextension}, using convex cocompactness, one verifies that this sub-bundle satisfies the uniform
flaring condition, and so is uniformly hyperbolic.  Moreover, it is also shown to be uniformly quasiisometrically
embedded in $Y$.  A similar construction can be done for three or more points on $S^1_\infty$, giving
ideal-polygon--sub-bundles of $Y$ (which are quasiisometric to tree bundles) with the same result.

Hyperbolicity of $Y$, and hence $\Gamma_G$, is arrived at by finding a sufficiently rich family of quasigeodesics connecting every pair of points, as is done for example in \cite{MM1} and \cite{hamenstadtcc}.   Given a pair of points $x,y \in Y$, to find the quasigeodesic in this family connecting $x$ to $y$, one picks the geodesic in an appropriate ideal-triangle--sub-bundle of $Y$ containing $x$ and $y$.  That this family is sufficiently rich is proved by appealing to the quasiisometric embedding and hyperbolicity properties of the sub-bundles as described above.
\end{proof}


If $F$ is a nonabelian Schottky group, then Farb and Mosher have proven \cite{FMII} that the extension $\Gamma_{\! F}$
is extremely quasiisometrically rigid---when $F$ is cyclic, the extension is quasiisometric to $\mathbb H^3$ and as
there are incommensurable hyperbolic $3$--manifolds that fiber over the circle with fiber $S$, we can hope for no such
theorem in that case. They prove

\begin{theorem}[Farb--Mosher \cite{FMII}] Let $F$ be a nonabelian Schottky group and let $H$ be a finitely generated group quasiisometric to $\Gamma_{\! F}$.  Then there is a finite normal subgroup $K$ of $H$ such that $H/K$ is commensurable to $\Gamma_{\! F}$.\qed
\end{theorem}

And the strong

\begin{theorem}[Farb--Mosher \cite{FMII}] Let $F$ be nonabelian and $Schottky$. The natural map from the abstract commensurator $\mathrm{Comm}(\Gamma_{\! F})$ of $\Gamma_{\! F}$ to the group $\mathrm{QI}(\Gamma_{\! F})$ of quasiisometries of $\Gamma_{\! F}$ is an isomorphism, and $\Gamma_{\! F}$ has finite index in $\mathrm{QI}(\Gamma_{\! F})$. \qed
\end{theorem}

\section{Gromov's ``hyperbolization'' question}

Thurston's Geometrization Conjecture predicts, in particular, that a closed aspherical $3$--manifold admits a
hyperbolic structure if and only if its fundamental group does not contain a copy of $\mathbb Z \oplus \mathbb Z$ (G.
Perelman has announced a proof of this conjecture \cite{perelman1,perelman2,perelman3}). There is a coarse version of
this conjecture for finitely generated groups, which we now describe.

For natural numbers $p$ and $q$, the \textit{Baumslag--Solitar group} $\textrm{BS}(p,q)$ is given by the presentation
\[
\textrm{BS}(p,q) = \langle a, b\ |\ a^{-1}b^pa=b^q \rangle.
\]
For geometric reasons, a hyperbolic group cannot contain any $\textrm{BS}(p,q)$. A question due to M. Gromov, arises, see \cite{bestvinalist}: is a finitely presented group with no Baumslag--Solitar subgroups hyperbolic? N. Brady has answered this question in the negative \cite{brady}. The counterexamples are subgroups of hyperbolic groups (and hence contain no $\textrm{BS}(p,q)$) and possess no finite $\mathrm K(\Gamma,1)$ (and hence fail to be hyperbolic as an appropriate Rips complex would serve as $\mathrm K(\Gamma,1)$). The refinement is

\begin{gromov}\label{gromov} If there is a finite $\mathrm K(\Gamma,1)$  and $\Gamma$ contains no Baumslag--Solitar subgroups, is $\Gamma$ hyperbolic?
\end{gromov}
\noindent The $\mathrm K(\Gamma,1)$ hypothesis is the natural analog of the aspherical hypothesis for \linebreak $3$--manifolds. Both hypotheses are satisfied by a great number of surface group extensions, as the following theorem illustrates.

\begin{theorem}\label{purelypA} Let $G$ be a subgroup of $\Mod(S)$. If $G$ is purely pseudo-Anosov, then $\Gamma_G$ contains no Baumslag--Solitar subgroups. If there is a finite $\mathrm{K}(G,1)$, then there is a finite $\mathrm{K}(\Gamma_G,1)$.
\end{theorem}

We begin with a lemma.

\begin{lemma} \label{bs} If $\mathrm{BS}(p,q)$ embeds in $\Mod(S)$, then $p=q$.
\end{lemma}
\begin{proof}
Since $b$ is an infinite order element of $\Mod(S)$, there is a pure power $b^k$ so that the Thurston decomposition of
$b^k$ is a composition of Dehn twists and pseudo-Anosov homeomorphisms all supported on pairwise disjoint connected
subsurfaces; see \cite{ivanov}. The powers of Dehn twists and the dilatations of pseudo-Anosovs are well defined and
invariant under conjugacy, and since $a^{-1}b^{pk}a = b^{qk}$, it follows that $p = q$.
\end{proof}

We note that there are nonabelian examples of groups $\mathrm{BS}(p,p)$ in $\Mod(S)$.  For example, let $\beta$ be an
essential simple closed curve, let $h$ be an element of order $p$ supported on $S - \beta$, and let $f$ be a partial
pseudo-Anosov supported on $S - \beta$ such that the subgroup $\langle h,f \rangle$ is isomorphic to $\mathbb Z/ \!
p\mathbb Z * \mathbb Z$, which we may do by Min's theorem (Theorem \ref{min} here). Now, let $a = f$, $b= hT_\beta$ and
notice that $a^{-1}b^pa=b^p$.   Considering normal forms, one can verify that the HNN--extension $\mathrm{BS}(p,p)$
embeds.


\begin{proof}[Proof of Theorem \ref{purelypA}]
Suppose that every non-identity element of $G$ is pseudo-Anosov and that $\Gamma_G$ contains the Baumslag--Solitar
group $\mathrm{BS}(p,q)$.  The group $\Gamma_G$ is a subgroup of $\Mod(\dot S)$, and so $p=q$ by Lemma \ref{bs}.  The
subgroup of $\mathrm{BS}(p,p)$ generated by $a$ and $b^p$ is isomorphic to $\mathbb Z \oplus \mathbb Z$ and its image
in $G$ must be nontrivial, as a hyperbolic surface group contains no $\mathbb Z \oplus \mathbb Z$ subgroups.
Furthermore, the image must be infinite cyclic as the centralizer of a pseudo-Anosov mapping class is
virtually cyclic \cite{mccarthythesis} and $G$ is torsion free.  We conclude that there is a pseudo-Anosov mapping class preserving the
conjugacy class of a nontrivial element of $\pi_1(S)$, which is impossible.

Let $\M_G$ be the cover of the moduli space $\M(S)$ corresponding to $G$. Supposing there exists a finite
$\mathrm{K}(G,1)$, the group $G$ is torsion free.
So, there is a canonical $S$--bundle
\begin{equation} \label{quo2}
S \to \mathcal B_G \to \M_G
\end{equation}
over $\M_G$ corresponding to the extension
\begin{equation} \label{ses2}
1 \to \pi_1(S) \to \Gamma_G \to G \to 1.
\end{equation}

More precisely, by a result of Bers \cite{bersfiber}, $\T(\dot S)$ fibers over $\T(S)$ with fiber the universal cover $\mathbb H^2$ of $S$:
\begin{equation} \label{uni1}
\mathbb H^2 \to \T(\dot S) \to \T(S)
\end{equation}
such that $\Mod(\dot S)$ acts on $\T(\dot S)$ by bundle maps.
This $\Mod(\dot S)$--action descends to a $\Mod(\dot S)$--action on $\T(S)$, and the descent factors through the natural action of $\Mod(S)$ via
the projection in (\ref{ses1}).  Furthermore, $\pi_1(S)$ acts on the fibers as the group of deck transformations of the universal cover $\mathbb H^2 \to S$.  The
universal curve (\ref{quo1}) is the quotient of (\ref{uni1}) by (\ref{ses1}), and (\ref{quo2}) is the quotient of (\ref{uni1}) by (\ref{ses2}).

Since $G$ is torsion free and the action on $\T(S)$ is properly discontinuous it follows that $G$ acts freely and so
$\M_G$ is a manifold. Teichm\"uller's Theorem asserts that $\T(S)$ is homeomorphic to a cell, and hence $\M_G$ is
aspherical.

We thus obtain a map $\mathrm{K}(G,1) \to \M_G$ and a pullback bundle
\[
S \to \mathcal K_G \to \mathrm{K}(G,1)
\]
with compact total space. The homotopy long exact sequence of the bundle reveals that $\mathcal K_G$ is aspherical.

We may choose a finite covering of $\mathcal K_G$ by absolute retracts whose nerve is homotopy equivalent to $\mathcal K_G$, by Weil's theorem (see pages 466--468 of \cite{corners}), and so obtain a finite $\mathrm K(\Gamma_G,1)$.
\end{proof}

As mentioned above, the only known purely pseudo-Anosov subgroups of $\Mod(S)$ are free.  Our final theorem shows that if there were purely pseudo-Anosov subgroups isomorphic to certain extensions, we would not have to check for convex cocompactness in order to answer Gromov's question in the negative.

\begin{theorem} Let $1 \to A \to B \to C \to 1$ be a short exact sequence of infinite groups such that $A$ and $B$ possess finite Eilenberg--Mac Lane spaces.  If $B$ is isomorphic to a purely pseudo-Anosov subgroup of $\Mod(S)$ (for $S$ closed), then $\Gamma_A$ and $\Gamma_B$ contain no Baumslag--Solitar subgroups, possess finite Eilenberg--Mac Lane spaces, and one of $\Gamma_A$ or $\Gamma_B$ is not hyperbolic.
\end{theorem}
\begin{proof} The first two statements follow from Theorem \ref{purelypA}.  If $B$ is not convex cocompact, then $\Gamma_B$ is not hyperbolic by Theorem \ref{ham}.  If $B$ is convex cocompact, then the weak hulls $\WH_A$ and $\WH_B$ are equal, as $A$ is normal in $B$.  Since $C$ is infinite, the index of $A$ in $B$ is infinite, and so $A$ could not act cocompactly on $\WH_B = \WH_A$.  By Theorem \ref{cocompactonhull}, $A$ is not convex cocompact, and so $\Gamma_A$ fails to be hyperbolic.
\end{proof}

\subsubsection*{Acknowledgements} The authors thank the Mathematics Department at the University of Michigan for its hospitality during the Colloquium, and Pete Storm for an enjoyable conversation that led to the discussion in Section \ref{measureconjsection}.  They would also like to thank Howard Masur and the referee for their interest and for carefully reading the paper and making many helpful suggestions.

\bibliographystyle{plain}
\bibliography{modsurvey}

\bigskip

\noindent Department of Mathematics, Brown University, Providence, RI 02912 \newline \noindent
\texttt{rkent@math.brown.edu}

\bigskip

\noindent Department of Mathematics, University of Illinois, Urbana-Champaign, IL 61801 \newline \noindent  \texttt{clein@math.uiuc.edu}

\end{document}

%% file: Ltable.pstex_t
\begin{picture}(0,0)%
\includegraphics{Ltable.pstex}%
\end{picture}%
\setlength{\unitlength}{4144sp}%
\begingroup\makeatletter\ifx\SetFigFont\undefined%
\gdef\SetFigFont#1#2#3#4#5{%
  \reset@font\fontsize{#1}{#2pt}%
  \fontfamily{#3}\fontseries{#4}\fontshape{#5}%
  \selectfont}%
\fi\endgroup%
\begin{picture}(2170,1671)(61,-915)
\put(571,609){\makebox(0,0)[lb]{\smash{{\SetFigFont{12}{14.4}{\rmdefault}{\mddefault}{\updefault}{\color[rgb]{0,0,0}$1$}%
}}}}
\put( 76,-106){\makebox(0,0)[lb]{\smash{{\SetFigFont{12}{14.4}{\rmdefault}{\mddefault}{\updefault}{\color[rgb]{0,0,0}$a$}%
}}}}
\put(2216,-351){\makebox(0,0)[lb]{\smash{{\SetFigFont{12}{14.4}{\rmdefault}{\mddefault}{\updefault}{\color[rgb]{0,0,0}$1$}%
}}}}
\put(1176,-851){\makebox(0,0)[lb]{\smash{{\SetFigFont{12}{14.4}{\rmdefault}{\mddefault}{\updefault}{\color[rgb]{0,0,0}$b$}%
}}}}
\end{picture}%

%% file: modsurveyfin.bbl
\def\cprime{$'$}
\begin{thebibliography}{10}

\bibitem{agol}
Ian Agol.
\newblock {Tameness of hyperbolic 3-manifolds, Preprint}.
\newblock \texttt{arXiv:math.GT/0405568},.

\bibitem{ahlforsmeasure}
Lars~V. Ahlfors.
\newblock Finitely generated {K}leinian groups.
\newblock {\em Amer. J. Math.}, 86:413--429, 1964.

\bibitem{ADPlconnect}
R.~C. Alperin, Warren Dicks, and J.~Porti.
\newblock The boundary of the {G}ieseking tree in hyperbolic three-space.
\newblock {\em Topology Appl.}, 93(3):219--259, 1999.

\bibitem{aas}
James~W. Anderson, Javier Aramayona, and Kenneth~J. Shackleton.
\newblock {A simple criterion for non-relative hyperbolicity and one-endedness
  of groups, {P}reprint}.
\newblock \texttt{arXiv:math.GT/0504271}.

\bibitem{AMlconnect}
James~W. Anderson and Bernard Maskit.
\newblock On the local connectivity of limit set of {K}leinian groups.
\newblock {\em Complex Variables Theory Appl.}, 31(2):177--183, 1996.

\bibitem{BeMa}
Alan~F. Beardon and Bernard Maskit.
\newblock Limit points of {K}leinian groups and finite sided fundamental
  polyhedra.
\newblock {\em Acta Math.}, 132:1--12, 1974.

\bibitem{bdm}
Jason Behrstock, Cornelia Drutu, and Lee Mosher.
\newblock {Thick metric spaces, relative hyperbolicity, and quasi-isometric
  rigidity, {P}reprint}.
\newblock \texttt{arXiv:math.GT/0512592}.

\bibitem{behrstockminsky}
Jason~A. Behrstock and Yair~N. Minsky.
\newblock {Dimension and rank for mapping class groups, {P}reprint}.
\newblock \texttt{arXiv:math.GT/0512352}.

\bibitem{bersfiber}
Lipman Bers.
\newblock Fiber spaces over {T}eichm\"uller spaces.
\newblock {\em Acta. Math.}, 130:89--126, 1973.

\bibitem{bestvinalist}
M.~Bestvina.
\newblock Questions in geometric group theory.
\newblock Available at \texttt{http://www.math.utah.edu/$\sim$bestvina/}.

\bibitem{combinationtheorem}
Mladen Bestvina and Mark Feighn.
\newblock A combination theorem for negatively curved groups.
\newblock {\em J. Differential Geom.}, 35(1):85--101, 1992.

\bibitem{birman}
Joan~S. Birman.
\newblock Mapping class groups and their relationship to braid groups.
\newblock {\em Comm. Pure Appl. Math.}, 22:213--238, 1969.

\bibitem{corners}
A.~Borel and J.-P. Serre.
\newblock Corners and arithmetic groups.
\newblock {\em Comment. Math. Helv.}, 48:436--491, 1973.
\newblock Avec un appendice: Arrondissement des vari\'et\'es \`a coins, par A.
  Douady et L. H\'erault.

\bibitem{bowditchstacks}
B.~Bowditch.
\newblock Stacks of hyperbolic spaces and ends of $3$--manifolds.
\newblock Preprint. Available at
  \texttt{www.maths.soton.ac.uk/staff/Bowditch/preprints.html}.

\bibitem{bowditchtop}
Brian~H. Bowditch.
\newblock A topological characterisation of hyperbolic groups.
\newblock {\em J. Amer. Math. Soc.}, 11(3):643--667, 1998.

\bibitem{bowditchcomplexofcurves}
Brian~H. Bowditch.
\newblock Hyperbolic 3-manifolds and the geometry of the curve complex.
\newblock In {\em European Congress of Mathematics}, pages 103--115. Eur. Math.
  Soc., Z\"urich, 2005.

\bibitem{brady}
Noel Brady.
\newblock Branched coverings of cubical complexes and subgroups of hyperbolic
  groups.
\newblock {\em J. London Math. Soc. (2)}, 60(2):461--480, 1999.

\bibitem{br}
Br. Brahmachaitanya.
\newblock {Cannon-Thurston Maps and Kleinian Groups: Amalgamation Geometry and
  the 5-holed Sphere, {P}reprint }.
\newblock \texttt{arXiv:math.GT/0512539}.

\bibitem{br1}
Br. Brahmachaitanya.
\newblock {Cannon-Thurston Maps and Kleinian Groups I: Pared Manifolds of
  Bounded Geometry, {P}reprint }.
\newblock \texttt{arXiv:math.GT/0503581 }.

\bibitem{br2}
Br. Brahmachaitanya.
\newblock {Cannon-Thurston Maps and Kleinian Groups II: i-bounded Geometry and
  a theorem of McMullen, {P}reprint }.
\newblock \texttt{arXiv:math.GT/0511104}.

\bibitem{brockfarb}
Jeffrey Brock and Benson Farb.
\newblock Curvature and rank of {T}eichm\"uller space.
\newblock {\em Amer. J. Math.}, 128(1):1--22, 2006.

\bibitem{calgabai}
Danny Calegari and David Gabai.
\newblock Shrinkwrapping and the taming of hyperbolic 3-manifolds.
\newblock {\em J. Amer. Math. Soc.}, 19(2):385--446 (electronic), 2006.

\bibitem{calta}
Kariane Calta.
\newblock Veech surfaces and complete periodicity in genus two.
\newblock {\em J. Amer. Math. Soc.}, 17(4):871--908 (electronic), 2004.

\bibitem{canaryends}
Richard~D. Canary.
\newblock Ends of hyperbolic {$3$}-manifolds.
\newblock {\em J. Amer. Math. Soc.}, 6(1):1--35, 1993.

\bibitem{duchin}
Moon Duchin.
\newblock Thin triangles and a multiplicative ergodic theorem for
  {T}eichm{\"u}ller geometry.
\newblock Preprint, \texttt{arXiv:math.GT/0508046}.

\bibitem{relativefarb}
B.~Farb.
\newblock Relatively hyperbolic groups.
\newblock {\em Geom. Funct. Anal.}, 8(5):810--840, 1998.

\bibitem{FMcc}
Benson Farb and Lee Mosher.
\newblock Convex cocompact subgroups of mapping class groups.
\newblock {\em Geom. Topol.}, 6:91--152 (electronic), 2002.

\bibitem{FMII}
Benson Farb and Lee Mosher.
\newblock The geometry of surface-by-free groups.
\newblock {\em Geom. Funct. Anal.}, 12(5):915--963, 2002.

\bibitem{floydlconnect}
William~J. Floyd.
\newblock Group completions and limit sets of {K}leinian groups.
\newblock {\em Invent. Math.}, 57(3):205--218, 1980.

\bibitem{gersten}
S.~M. Gersten.
\newblock Cohomological lower bounds for isoperimetric functions on groups.
\newblock {\em Topology}, 37(5):1031--1072, 1998.

\bibitem{hamenstadtcc}
Ursula Hamenst{\"a}dt.
\newblock Geometry of the complex of curves and of {T}eichm{\"u}ller space,
  {P}reprint.
\newblock \texttt{arXiv:math.GT/0502256}.

\bibitem{hamenstadt}
Ursula Hamenst{\"a}dt.
\newblock {Word hyperbolic extensions of surface groups}.
\newblock Preprint, \texttt{arXiv:math.GT/0505244}.

\bibitem{hamenstadtrank}
Ursula Hamenstaedt.
\newblock {Geometry of the mapping class groups III: Geometric rank,
  {P}reprint}.
\newblock \texttt{arXiv:math.GT/0512429}.

\bibitem{hubbardmasur}
John Hubbard and Howard Masur.
\newblock Quadratic differentials and foliations.
\newblock {\em Acta Math.}, 142(3-4):221--274, 1979.

\bibitem{KL}
Richard P.~Kent IV and Christopher~J. Leininger.
\newblock {Shadows of mapping class groups: capturing convex co-compactness,
  {P}reprint}.
\newblock \texttt{arXiv:math.GT/0505114}.

\bibitem{KLtrip}
Richard P.~Kent IV and Christopher~J. Leininger.
\newblock {Uniform convergence in the mapping class group}.
\newblock In progress.

\bibitem{ivanov}
Nikolai~V. Ivanov.
\newblock {\em Subgroups of {T}eichm\"uller modular groups}, volume 115 of {\em
  Translations of Mathematical Monographs}.
\newblock American Mathematical Society, Providence, RI, 1992.
\newblock Translated from the Russian by E. J. F. Primrose and revised by the
  author.

\bibitem{kn}
Anders Karlsson and Guennadi~A. Noskov.
\newblock Some groups having only elementary actions on metric spaces with
  hyperbolic boundaries.
\newblock {\em Geom. Dedicata}, 104:119--137, 2004.

\bibitem{kerckhoffhandle}
Steven~P. Kerckhoff.
\newblock The measure of the limit set of the handlebody group.
\newblock {\em Topology}, 29(1):27--40, 1990.

\bibitem{klarreichlconnect}
Erica Klarreich.
\newblock Semiconjugacies between {K}leinian group actions on the {R}iemann
  sphere.
\newblock {\em Amer. J. Math.}, 121(5):1031--1078, 1999.

\bibitem{weilletter}
Martin~H. Krieger.
\newblock A 1940 letter of {A}ndr\'e {W}eil on analogy in mathematics.
\newblock {\em Notices Amer. Math. Soc.}, 52(3):334--341, 2005.
\newblock Excerpted from {\it Doing mathematics} [World Scientific Publishing
  Co., Inc., River Edge, NJ, 2003; MR1961400].

\bibitem{LRcombination}
C.~J. Leininger and A.~W. Reid.
\newblock A combination theorem for {V}eech subgroups of the mapping class
  group.
\newblock {\em Geom. Funct. Anal.}, 16(2):403--436, 2006.

\bibitem{graphsofveech}
Christopher~J. Leininger.
\newblock Graphs of {V}eech groups and ideal boundaries.
\newblock Work in progress.

\bibitem{leiningertwist}
Christopher~J. Leininger.
\newblock On groups generated by two positive multi-twists: {T}eichm\"uller
  curves and {L}ehmer's number.
\newblock {\em Geom. Topol.}, 8:1301--1359 (electronic), 2004.

\bibitem{masurclass}
Howard Masur.
\newblock On a class of geodesics in {T}eichm\"uller space.
\newblock {\em Ann. of Math. (2)}, 102(2):205--221, 1975.

\bibitem{masuruniquely}
Howard Masur.
\newblock Uniquely ergodic quadratic differentials.
\newblock {\em Comment. Math. Helv.}, 55(2):255--266, 1980.

\bibitem{masurinterval}
Howard Masur.
\newblock Interval exchange transformations and measured foliations.
\newblock {\em Ann. of Math. (2)}, 115(1):169--200, 1982.

\bibitem{twoboundaries}
Howard Masur.
\newblock Two boundaries of {T}eichm\"uller space.
\newblock {\em Duke Math. J.}, 49(1):183--190, 1982.

\bibitem{masurhandle}
Howard Masur.
\newblock Measured foliations and handlebodies.
\newblock {\em Ergodic Theory Dynam. Systems}, 6(1):99--116, 1986.

\bibitem{masurhaus}
Howard Masur.
\newblock Hausdorff dimension of the set of nonergodic foliations of a
  quadratic differential.
\newblock {\em Duke Math. J.}, 66(3):387--442, 1992.

\bibitem{masurtabachnikov}
Howard Masur and Serge Tabachnikov.
\newblock Rational billiards and flat structures.
\newblock In {\em Handbook of dynamical systems, Vol.\ 1A}, pages 1015--1089.
  North-Holland, Amsterdam, 2002.

\bibitem{MM1}
Howard~A. Masur and Yair~N. Minsky.
\newblock Geometry of the complex of curves. {I}. {H}yperbolicity.
\newblock {\em Invent. Math.}, 138(1):103--149, 1999.

\bibitem{MM2}
Howard~A. Masur and Yair~N. Minsky.
\newblock Geometry of the complex of curves. {II}. {H}ierarchical structure.
\newblock {\em Geom. Funct. Anal.}, 10(4):902--974, 2000.

\bibitem{MMunstable}
Howard~A. Masur and Yair~N. Minsky.
\newblock Unstable quasi-geodesics in {T}eichm\"uller space.
\newblock In {\em In the tradition of Ahlfors and Bers (Stony Brook, NY,
  1998)}, volume 256 of {\em Contemp. Math.}, pages 239--241. Amer. Math. Soc.,
  Providence, RI, 2000.

\bibitem{MMdiskset}
Howard~A. Masur and Yair~N. Minsky.
\newblock Quasiconvexity in the curve complex.
\newblock In {\em In the tradition of Ahlfors and Bers, III}, volume 355 of
  {\em Contemp. Math.}, pages 309--320. Amer. Math. Soc., Providence, RI, 2004.

\bibitem{masurwolf}
Howard~A. Masur and Michael Wolf.
\newblock Teichm\"uller space is not {G}romov hyperbolic.
\newblock {\em Ann. Acad. Sci. Fenn. Ser. A I Math.}, 20(2):259--267, 1995.

\bibitem{mctits}
John McCarthy.
\newblock A ``{T}its-alternative'' for subgroups of surface mapping class
  groups.
\newblock {\em Trans. Amer. Math. Soc.}, 291(2):583--612, 1985.

\bibitem{mccarthypapa}
John McCarthy and Athanase Papadopoulos.
\newblock Dynamics on {T}hurston's sphere of projective measured foliations.
\newblock {\em Comment. Math. Helv.}, 64(1):133--166, 1989.

\bibitem{mccarthythesis}
John~D. McCarthy.
\newblock Normalizers and centralizers of pseudo-anosov mapping classes.
\newblock Columbia University Ph.D. thesis.

\bibitem{mcmullenlconnect}
Curtis~T. McMullen.
\newblock Local connectivity, {K}leinian groups and geodesics on the blowup of
  the torus.
\newblock {\em Invent. Math.}, 146(1):35--91, 2001.

\bibitem{mcmullenbilliard}
Curtis~T. McMullen.
\newblock Billiards and {T}eichm\"uller curves on {H}ilbert modular surfaces.
\newblock {\em J. Amer. Math. Soc.}, 16(4):857--885 (electronic), 2003.

\bibitem{min}
Honglin Min.
\newblock Work in progress.

\bibitem{minskylconnect}
Yair~N. Minsky.
\newblock On rigidity, limit sets, and end invariants of hyperbolic
  {$3$}-manifolds.
\newblock {\em J. Amer. Math. Soc.}, 7(3):539--588, 1994.

\bibitem{minskyextremal}
Yair~N. Minsky.
\newblock Extremal length estimates and product regions in {T}eichm\"uller
  space.
\newblock {\em Duke Math. J.}, 83(2):249--286, 1996.

\bibitem{minskycrelle}
Yair~N. Minsky.
\newblock Quasi-projections in {T}eichm\"uller space.
\newblock {\em J. Reine Angew. Math.}, 473:121--136, 1996.

\bibitem{minskyinj}
Yair~N. Minsky.
\newblock Kleinian groups and the complex of curves.
\newblock {\em Geom. Topol.}, 4:117--148, 2000.

\bibitem{minskybound}
Yair~N. Minsky.
\newblock Bounded geometry for {K}leinian groups.
\newblock {\em Invent. Math.}, 146(1):143--192, 2001.

\bibitem{mosherextension}
Lee Mosher.
\newblock Hyperbolic extensions of groups.
\newblock {\em J. Pure Appl. Algebra}, 110(3):305--314, 1996.

\bibitem{hypbyhyp}
Lee Mosher.
\newblock A hyperbolic-by-hyperbolic hyperbolic group.
\newblock {\em Proc. Amer. Math. Soc.}, 125(12):3447--3455, 1997.

\bibitem{mosherstable}
Lee Mosher.
\newblock Stable {T}eichm\"uller quasigeodesics and ending laminations.
\newblock {\em Geom. Topol.}, 7:33--90 (electronic), 2003.

\bibitem{moshersurvey}
Lee Mosher.
\newblock Problems in the geometry of surface group extensions.
\newblock In {\em Problems on mapping class groups and related topics},
  volume~74 of {\em Proc. Symp. Pure and Applied Math.}, pages 261--273. 2006.

\bibitem{perelman3}
Grisha Perelman.
\newblock {Finite extinction time for the solutions to the Ricci flow on
  certain three--manifolds, {P}erlman}.
\newblock \texttt{arXiv:math.DG/0307245}.

\bibitem{perelman2}
Grisha Perelman.
\newblock {Ricci flow with surgery on three--manifolds, {P}reprint}.
\newblock \texttt{arXiv:math.DG/0303109}.

\bibitem{perelman1}
Grisha Perelman.
\newblock {The entropy formula for the Ricci flow and its geometric
  applications, {P}reprint}.
\newblock \texttt{arXiv:math.DG/0211159}.

\bibitem{rafi}
Kasra Rafi.
\newblock A characterization of short curves of a {T}eichm\"uller geodesic.
\newblock {\em Geom. Topol.}, 9:179--202, 2005.

\bibitem{swenson}
Eric~L. Swenson.
\newblock Quasi-convex groups of isometries of negatively curved spaces.
\newblock {\em Topology Appl.}, 110(1):119--129, 2001.
\newblock Geometric topology and geometric group theory (Milwaukee, WI, 1997).

\bibitem{weilletterfrench}
Andr{\'e} Weil.
\newblock {\em Scientific works. {C}ollected papers. {V}ol. {I} (1926--1951)}.
\newblock Springer-Verlag, New York, 1979.

\end{thebibliography}
